# Prime Numbers, Dirichlet Density, and Benford's Law

Alex Ely Kossovsky*

**Abstract**

The Prime Numbers are well-known for their paradoxical stand regarding Benford's Law. On one hand they adamantly refuse to obey the law of Benford in the usual sense, namely that of a normal density of the proportion of primes with d as the leading digit, yet on the other hand, the Dirichlet density for the subset of all primes with d as the leading digit is indeed LOG(1 +1/d). In this article the superficiality of the Dirichlet density result is demonstrated and explained in terms of other well-known and established results in the discipline of Benford's Law, conceptually concluding that prime numbers cannot be considered Benford at all, in spite of the Dirichlet density result. In addition, a detailed examination of the digital behavior of prime numbers is outlined, showing a distinct digital development pattern, from a slight preference for low digits at the start for small primes, to a complete digital equality for large primes in the limit as the prime number sequence goes to infinity. Finally an exact analytical expression for the density of the logarithms of primes is derived and shown to be always on the rise at the macro level, an observation that is also confirmed empirically.

* Email: akossovsky@gmail.com

# SECTION 1: Empirical Checks on Distributions of Primes and Digits
-----------------------------------------------------------------------------------------------------------

Prime numbers have intrigued mathematicians for millennia. While the integers are considered to be relatively well-understood, not enough is known about the primes. Indeed there are still many open problems concerning the primes, and this renders them a bit mysterious. Much progress has been made regarding the prime number distribution since Gauss initial conjectures. Although the prime number distribution is deterministic in the sense that exact rules determine whether an integer is prime or composite, its apparent randomness has caused many to speculate about possible stochastic interpretations. On one hand, locally for small ranges of integers, prime numbers seem to be randomly distributed with no other law or pattern than that of chance. But on the other hand, globally for large ranges of integers, the distribution of primes tells of a remarkably predictable pattern. This tension between micro randomness and macro order renders the distribution of primes a fascinating problem indeed.

First Leading Digit or First Significant Digit is the first (non-zero) digit of a given number appearing on the left-most side. For 567.34 the leading digit is 5. For 0.0367 the leading digit is 3, as we discard the zeros. For the lone integer 6 the leading digit is 6. For negative numbers we simply discard the sign, hence for -62.97 the leading digit is 6.

Naturally, when digit d appears first in a number composed of several digits, we call d the 'leader', as if it leads all the other digits trailing behind it to the right.

**6**13 → digit 6
0.000**2**867 → digit 2
 **7** → digit 7
-**7** → digit 7
**1**,653,832 → digit 1
-0.**4**56398 → digit 4

Perhaps it is tempting to intuit that for numbers in typical real life data sets, all nine digits {1, 2, 3, 4, 5, 6, 7, 8, 9} should be equally likely to occur and thus uniformly distributed. However this is not the case. In fact, low digits such as {1, 2, 3} actually occur with very high frequencies within the 1st place position of typical everyday data, while high digits such as {7, 8, 9} have very little overall proportion of occurrence. So much so that the proportion of everyday typical numbers starting with digit 1 is about seven times the proportion of numbers starting with digit 9! About 30% of typical everyday numbers in use start with digit 1, while only about 4% start with digit 9.

That low digits occur in the first order much more often than high digits can be seen directly from the analytical expression:

Probability[1st digit is d] = $LOG_{10}(1 + 1/d)$

This is known as Benford's Law. The implied set of proportions is known as the logarithmic distribution. Figure 1 depicts the distribution.

| Digit | Probability |
|---|---|
| 1 | 30.1% |
| 2 | 17.6% |
| 3 | 12.5% |
| 4 | 9.7% |
| 5 | 7.9% |
| 6 | 6.7% |
| 7 | 5.8% |
| 8 | 5.1% |
| 9 | 4.6% |

**Figure 1:** Benford's Law

Remarkably, Benford's Law is found in a bewildering variety of real life data sets relating to physics, chemistry, astronomy, economics, finance, accounting, geology, biology, engineering, and governmental census data.

As an example, geological data on time (in units of seconds) between earthquakes is highly logarithmic, with only tiny deviations from LOG(1 + 1/d). Figure 2 depicts a very small sample of 40 values randomly taken from the dataset on all global earthquake occurrences in the year 2012. 1st digits are emphasized in bold font in Figure 3.

```
Digit Index:                           {  1,   2,   3,   4,  5, 6, 7, 8,  9 }
Digits Count, totaling 40 values:      { 15,   8,   6,   4,  4, 0, 2, 1,  0 }
Proportions of Digits, '%' sign omitted: {38,  20,  15,  10, 10, 0, 5, 3,  0 }
```

| 285.29  | 185.35  | 2579.80 | 27.11   |
|---------|---------|---------|---------|
| 5330.22 | 1504.49 | 1764.41 | 574.46  |
| 1722.16 | 815.06  | 3686.84 | 1501.61 |
| 494.17  | 362.48  | 1388.13 | 1817.27 |
| 3516.80 | 5049.66 | 2414.06 | 387.78  |
| 4385.23 | 2443.98 | 2204.12 | 1224.42 |
| 1965.46 | 3.61    | 1347.30 | 271.23  |
| 3247.99 | 753.80  | 1781.45 | 593.59  |
| 1482.64 | 1165.04 | 4647.39 | 1219.19 |
| 251.12  | 7345.52 | 1368.79 | 4112.13 |

**Figure 2:** Random Sample of 40 Earthquake Values

| **2**85.29  | **1**85.35  | **2**579.80 | **2**7.11   |
|-------------|-------------|-------------|-------------|
| **5**330.22 | **1**504.49 | **1**764.41 | **5**74.46  |
| **1**722.16 | **8**15.06  | **3**686.84 | **1**501.61 |
| **4**94.17  | **3**62.48  | **1**388.13 | **1**817.27 |
| **3**516.80 | **5**049.66 | **2**414.06 | **3**87.78  |
| **4**385.23 | **2**443.98 | **2**204.12 | **1**224.42 |
| **1**965.46 | **3**.61    | **1**347.30 | **2**71.23  |
| **3**247.99 | **7**53.80  | **1**781.45 | **5**93.59  |
| **1**482.64 | **1**165.04 | **4**647.39 | **1**219.19 |
| **2**51.12  | **7**345.52 | **1**368.79 | **4**112.13 |

**Figure 3:** First Digits of the Earthquake Sample

The Prime Numbers can be viewed as a dynamic sequence, advancing from 2 forward to higher values, in a conceptually similar sense to exponential growth series or random log walk; as opposed to viewing them as a static collection of unique numbers.

When one compares how the prime numbers sequence advances along the x-axis to exponential growth series, it immediately becomes apparent that it strongly resembles very slow exponential growth rates, rarely passing an **Integral Power Of Ten** (IPOT) point, and it keeps slowing down, growing ever more slowly, constantly shifting to lower growth rates. Such state of affairs suggests that any empirical trials regarding primes compliance or non-compliance with Benford's Law should all stop at IPOT points.
In other words, it is essential that we take the digital pulse of prime numbers on intervals such that an unbiased examination is performed where all possible first digits are equiprobable a priori. A superior approach regarding this boundary issue is not only to pay close attention to IPOT points in all empirical tests of finite subsets of prime numbers, but also to restrict such tests to intervals standing between two adjacent IPOT values, namely intervals on $(10^{INTEGER}, 10^{INTEGER + 1})$. To vividly demonstrate the importance of this approach in empirical trials, imagine a digital test performed on all the primes from 2 to 400, or on all the primes in the interval (100, 400). Clearly all such tests are biased, unfairly favoring digits {1, 2, 3}, while strongly discriminating against digits 4 to 9. The results of such misguided tests tell us more about the measuring rod itself (interval chosen) than about the measured object (digital configuration of primes). All such results are worthless as they do not convey any useful information.

Here are some proper empirical results:

| | |
|---|---|
| Primes between 2 and 100 | {16.0, 12.0, 12.0, 12.0, 12.0, 8.0, 16.0, 8.0, 4.0} |
| Primes between 2 and 1,000 | {14.9, 11.3, 11.3, 11.9, 10.1, 10.7, 10.7, 10.1, 8.9} |
| Primes between 2 and 10,000 | {13.0, 11.9, 11.3, 11.3, 10.7, 11.0, 10.2, 10.3, 10.3} |
| Primes between 2 and 100,000 | {12.4, 11.8, 11.4, 11.1, 11.0, 10.6, 10.7, 10.5, 10.5} |
| Primes between 2 and 1,000,000 | {12.2, 11.6, 11.4, 11.1, 11.0, 10.8, 10.7, 10.6, 10.5} |
| | |
| Primes between 10 and 100 | {19.0, 9.5, 9.5, 14.3, 9.5, 9.5, 14.3, 9.5, 4.8} |
| Primes between 100 and 1,000 | {14.7, 11.2, 11.2, 11.9, 9.8, 11.2, 9.8, 10.5, 9.8} |
| Primes between 1,000 and 10,000 | {12.7, 12.0, 11.3, 11.2, 10.7, 11.0, 10.1, 10.4, 10.6} |
| Primes between 10,000 and 100,000 | {12.4, 11.8, 11.5, 11.1, 11.0, 10.5, 10.8, 10.5, 10.5} |
| Primes between 100,000 & 1 million | {12.2, 11.6, 11.4, 11.1, 11.0, 10.8, 10.8, 10.6, 10.5} |

[**Note:** the vectors of nine elements refer to proportional percents, starting with digit 1 on the left, and ending with digit 9 on the right; the '%' sign is omitted for brevity.]

A clear picture emerges regarding digital behavior of the prime numbers. There is a consistent pattern here, slightly favoring low digits, but nothing like the dramatic digital skewness of Benford's Law is found here. It also appears as if there is a tendency towards a more equitable digital configuration as higher primes are considered.

For example, the proportions of digit 1 within adjacent IPOT intervals are constantly and monotonically decreasing: **19.0 → 14.7 → 12.7 → 12.4 → 12.2**. Does digit 1 aim at the ultimate digital equality of 11.1% (that is 1/9) in the limit as primes go to infinity? Unfortunately, the author has no access to a much bigger list of primes with values beyond one million integers. Hence, in order to learn more about the digital behavior of primes, only theoretical reasoning shall be employed here instead of further empirical examinations.

Figure 4 depicts the way the sequence of the 21 prime numbers from 11 to 97 marches along the log-axis. This is done by simply converting these 21 prime numbers into their (decimal) logarithm values and plotting them. The boundary of 11 and 97 is purposely selected so as to stand as close as possible between the IPOT postmarks of 10 and 100.

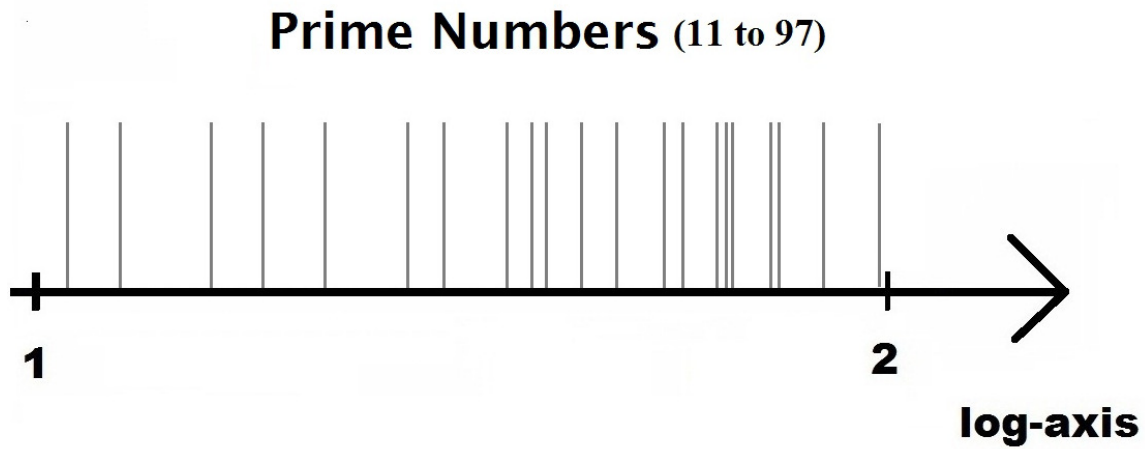

**Figure 4:** Decelerated March Along the log-axis – Primes 11 to 97

If one cannot visually ascertain from Figure 4 that this sequence slows down, becoming more and more concentrated towards the right, then the next example of Figure 5 would convince anyone of this generic pattern in all plots of logarithms of prime numbers. Figure 5 depicts the way the sequence of 143 prime numbers - beginning with 101 and ending with 997 - marches along the log-axis. A more detailed examination of the log values of these 143 primes is given in Figure 6 which focuses separately on the first half and on last half of the entire interval of (2.0, 3.0) for better visualization.

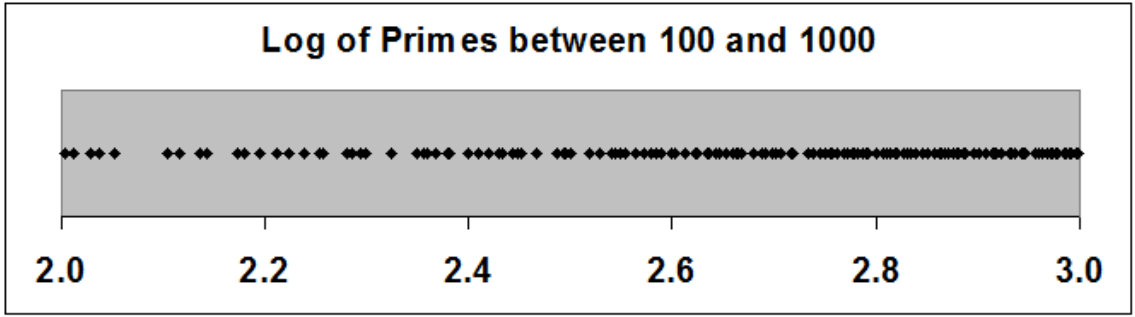

**Figure 5:** Decelerated March Along the log-axis – Primes 100 - 1000

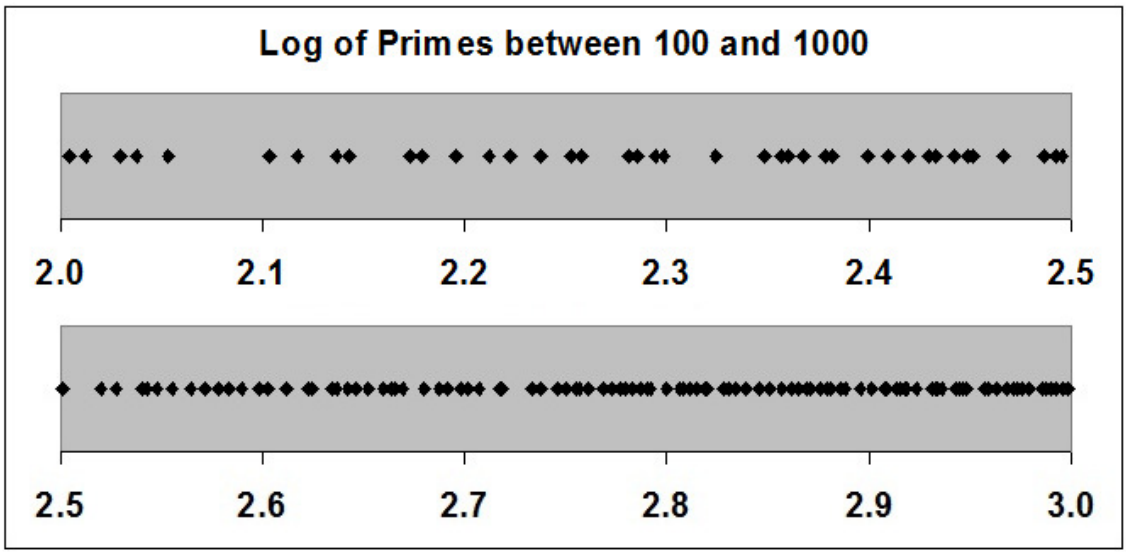

**Figure 6:** Decelerated log-axis March – Primes 100 - 1000 - Detailed

Clearly, distances on the log-axis between consecutive log values of primes generally decrease, and overall points are getting increasing more dense. The somewhat 'random' nature of how primes are spread among the integers implies that locally in some less frequent cases distances increase, but globally and most frequently they decrease. The implication of such log pattern is that the density of log of primes is rising on the right. Hence as discussed in Kossovsky (2014) chapter 63, and in particularly as shown in Figure 4.21 there, the expectation here is of an approximate digital equality, or at least of a much milder skewness in favoring low digits – just as was indeed observed and confirmed above empirically. This latest agreement and harmony between the generic theoretical understanding of the entire Benford phenomenon and the particular case of prime numbers (on top of numerous other cases) is highly reassuring.

Figures 7 and 8 pertain only to the particular interval standing between the adjacent IPOT values 10,000 and 100,000, and one <u>cannot</u> assume that these two figures demonstrate the generic way how all prime numbers behave between any consecutive IPOT points. A probable tendency towards digital equality is a distinct possibility as discussed earlier, and this issue has to be further investigated. At least on this adjacent IPOT interval, the histogram of the sequence of the primes themselves falls off gently and gradually to the right, while the histogram of the log of primes almost consistently rises.

There are 30 bins being utilized in the construction of the histogram of the prime numbers here, with each bin spanning (100,000-10,000)/30, or 90,000/30, namely the width of 3000 integers. Such global or aggregate perspective on how primes occur is almost smooth, yet it still shows some slight random fluctuations. A much sharper and more refined lens with a narrower bin width of only 50 for example, focusing on how many primes occur every 50 integers, would show many more random and wilder fluctuations in its histogram than the gentle ones seen in Figure 7. To mere mortals like humans, the fluctuations in the occurrences of primes do appear somewhat random and chaotic, but to the Goddess of Arithmetic all this appears rational, predictable, and even orderly, and she has proclaimed and determined all this long before the Big Bang noisily erupted up there in the sky 13.8 billion years ago or so.

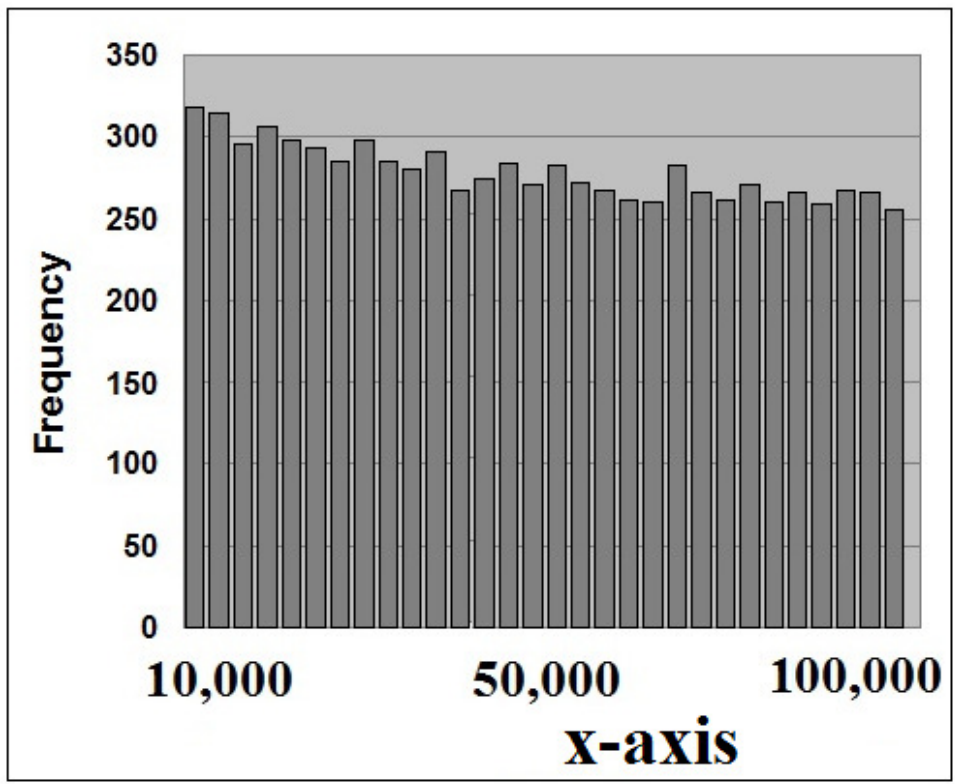

**Figure 7:** Histogram of Primes between 10,000 and 100,000

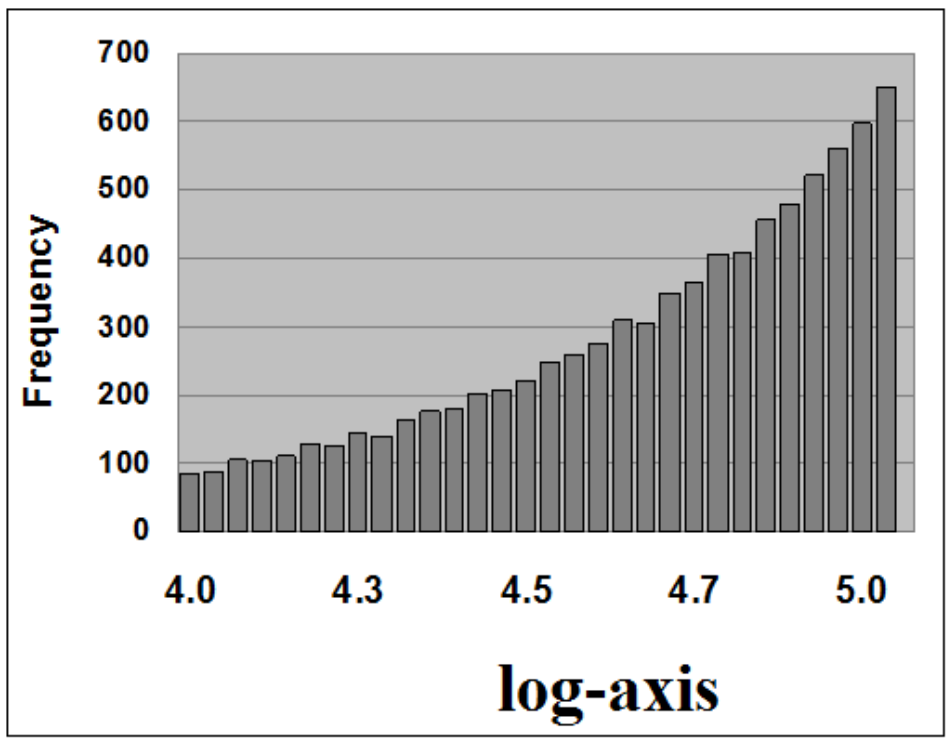

**Figure 8:** Histogram of Log of Primes between 10,000 and 100,000

As discussed in extreme generality in Kossovsky (2014) at the end of chapter 99, sequences that march along the log-axis in a coordinated, organized, and deliberate way with respect to the integer postmarks, exhibit digital behavior that may not be Benford (unless tiny steps are always taken covering almost all corners and sections of log-axis), while a random and chaotic log march almost always leads to Benford. This statement is merely the grand conceptual guideline, and randomness in the march along the log-axis does not immediately imply that mantissa is uniform, hence mathematicians are required to provide a rigorous proof for each particular case, showing whether or not it is logarithmic. For the case of the sequence of prime numbers, even though at the micro level the march appears random, at the macro level a clear and decisive pattern is seen of an overall rise along the log-axis between the (global) integer postmarks, and this pattern guarantees that mantissa is not at all uniform. As shall be demonstrated later, log of primes always rises, in the beginning as well as at infinity; hence Figure 8 depicts the true generic shape of all the histograms of log of primes, in contrast to Figure 7 of the primes themselves which evolves eventually into flat and uniform histograms at infinity. The huge number of primes residing between each pair of consecutive integer-postmarks on the log-axis renders those postmarks definitely global, affecting issues on the macro level. As the sequence of log of the primes passes more and more integer postmarks on its long march to infinity, mantissa is distributed in a rather skewed manner on its defined range of (0, 1), rising from 0 all the way to 1, depriving low digits of their usual strong Benfordian advantage, and granting high digits (almost and eventually complete) equitable proportions.

# SECTION 2:  Prime Number Theorem Hinting at Digital Configuration
-------------------------------------------------------------------------------------------------------------

Let us now rigorously prove that the leading digits of prime numbers are equally
distributed in the limit as they approach infinity, and also demonstrate that the generic
density of log of prime numbers is continuously rising throughout on the macro level,
from log(2) all the way to log(infinity). We shall establish both results by deducing them
from the asymptotic law of distribution of prime numbers, also widely known as
'**The prime Number Theorem**". The function π(x) is defined as the number of primes
less than or equal to x. Carl Friedrich Gauss is generally credited with first conjecturing
that π(x) is asymptotically x/ln(x). The notation ln(x) refers to the natural logarithm of x
with base e. The notation log(x) refers to the decimal logarithm of x with base 10. Gauss
arrived at the conjecture after receiving a book of tables of logarithms as a present at the
age of about 15 years in 1792. The book included as an appendix a table of primes,
intended just as a mathematical curiosity, and Gauss was able to make this unexpected
connection between primes and logarithms. Later in his life Gauss refined his estimate to
∫ 1/ln(t)dt with limits of integration running from 2 to x.

$$\pi(x) = Number\ of\ Primes\ up\ to\ x$$

$$\pi(x) \approx \int_{2}^{x} \frac{1}{\ln(t)}\ dt$$

This integral is strongly suggestive of the notion that the 'density' of primes in the
approximate vicinity around t should be 1/ln(t). By 'density' we mean the proportion of
consecutive integers centered around t that are primes. This is <u>not</u> the classic definition of
density in mathematical statistics where total area sums to 1, and where (dx)*(density) at
x indicates the (finite) proportion of all values falling within the tiny interval from x to
x + dx. This is the case for the prime numbers since there are infinitely many of them,
and one cannot state that a certain proportion of all the primes falls within any finite
x-axis subinterval. A better perspective about 1/ln(x) is to think of it as the 'Relative
Density', comparing counts of primes for a variety of ranges on the x-axis.

Prime Density = Proportion of primes within a range of consecutive integers
Prime Density = (Primes) / (Composites + Primes)
Prime Density = Prime Count / Integer Range
Prime Count =  Integer Range * Prime Density

Around the integer 373, the supposed 'density' of primes is 1/ln(373), or 0.169, meaning
that about **16.9%** of the integers around 373 are expected to be primes. Indeed, the
primes in this vicinity are {349, **353, 359, 367, <u>373</u>, 379, 383, 389**, 397}. Hence from
around 351 to around 393, having an interval the size of (393 – 351) = 42, we have 7
primes, and that yields the ratio 7/42, namely **16.7%**. In other words, out of 42 integers
with the potential of being primes only 7 turned out to be so, namely only 16.7% of them.

Around the integer 1709, the supposed 'density' of primes is 1/ln(1709), or 0.134, meaning that about **13.4%** of the integers around 1709 are expected to be primes. Indeed, the primes in this vicinity are {1669, **1693, 1697, 1699, 1709, 1721, 1723, 1733**, 1741}. Hence from around 1683 to around 1737, having an interval the size of (1737 – 1683) = 54, we have 7 primes, and that yields the ratio 7/54, namely **13.0%**. In other words, out of 54 integers with the potential of being primes only 7 turned out to be so, namely only 13.0% of them.

The two examples above should be further refined by considering wider ranges around the integers, instead of the narrow focus of only about 40 to 50 consecutive integers. The wider the range, the less randomness and fluctuations are observed, and the asymptotic law of distribution of prime numbers is better applied. A good choice perhaps is about 100 consecutive integers, where a remarkable fit is obtained between the theoretical prime density of 1/ln(x), and the empirically observed density. Figure 9 depicts the rather good fit between the theoretical 1/ln(x) density and the empirical (actual) prime density for the crude and narrow range of only 27 consecutive integers, examining 548 primes, from the prime 47 to the prime 4079.

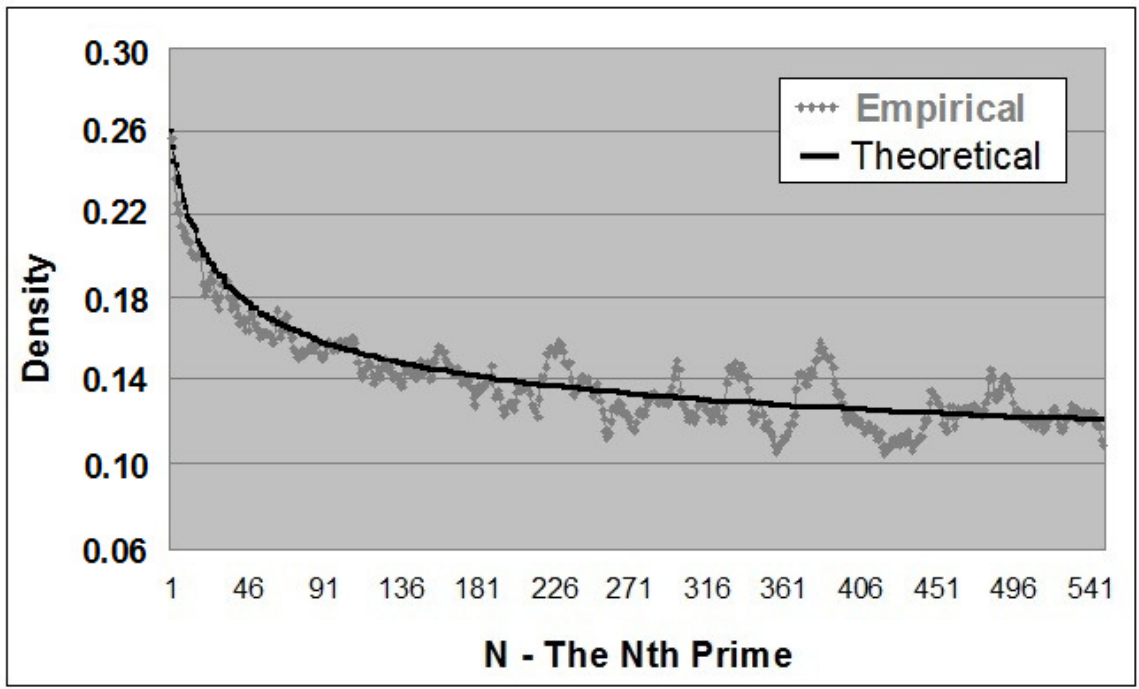

**Figure 9:** The Fit between Empirical Density and Theoretical 1/ln(x)

The above discussion neatly explains why the expression **x/ln(x)** conjectured by young Gauss for the total number of primes up to x consistently underestimates actual counts of primes; always falling short. This is so because the density of primes varies locally, falling as we consider higher integers, being inversely proportional to ln(x). The crude use of x/ln(x) estimates the entire number of integers from 2 to x as simply x, approximating (x - 2) as x, which is actually correct for large x. But it simple-mindedly and lazily estimates the density as a constant 1/ln(x) throughout, from 2 to x, although in reality the density at the last integer x is at its weakest point. Thus (x)*(1/ln(x)) is crudely expressing [Integer Range]*[Prime Density] for the value of [Prime Count]. Since density of primes for all integers less than x is a bit greater than 1/ln(x), such as in 1/ln(x – 568) or 1/ln(x – 37912), and so forth, this crude estimate is always less than the actual count of primes. Gauss the adult, considering the definitive integral, is willing to be diligent and break up the entire interval into many locals ones to account for the variability in the density and the fact that density is higher for lower integers, hence obtaining a much closer fit to the actual count of primes.

From the expression of the density of primes 1/ln(x) one can deduce that ultimately in the limit at infinity first digits are equality distributed within each interval standing between adjacent IPOT values, namely between **($10^N$, $10^{(N + 1)}$)**, where N is an integer.

Value of prime density at the left-most point:   $1/\ln(10^N)$
Value of prime density at the right-most point: $1/\ln(10^{(N + 1)})$

Ratio of left-most point density to right-most point density is:

$1/\ln(10^N) / 1/\ln(10^{(N + 1)}) = \ln(10^{(N + 1)}) / \ln(10^N)$

Transforming from base e to base 10 via $LOG_A(X) = LOG_B(X) / LOG_B(A)$:

$\ln(10^{(N + 1)}) / \ln(10^N) = \log(10^{(N + 1)})/\log(e) / \log(10^N)/\log(e) = \log(10^{(N + 1)}) / \log(10^N) =$
$(N+1)*\log(10) / N*\log(10) = (N+1)*1 / N*1 =$ **(N + 1)/N**

In the limit as primes go to infinity, N also goes to infinity albeit reluctantly and 'much more slowly' (N constitutes approximately the exponent of the largest primes on the way to infinity, not the values of the primes themselves) hence the ratio (N + 1)/N is 1 since the term 1 in the numerator is negligible in the grand scheme of things. In conclusion: leading digit distribution of prime numbers is asymptotically uniform on adjacent IPOT intervals, and therefore leading digit distribution of prime numbers is asymptotically uniform on the entire range to infinity.

Another way of looking at the density of primes is to place its gently falling curve within the series of the following falling curves in order of steepness:

$1/x^3$    $1/x^2$    $1/x$    $1/\text{square-root}(x)$    $1/\text{cube-root}(x)$    $1/\ln(x)$    $1/\text{constant}$

One can obtain all seven derivatives and make a comparison. Clearly $1/x^3$ is falling the steepest, the fastest, and 1/constant is falling the least; indeed, 1/constant is a uniform flat line, not falling at all. The threshold from which a definite integral from any real point k to infinity diverges is from $1/x$, and all the other curves on its right diverge as well, because in the limit their curves approach a flat line, implying an infinite area to infinity.

The author has access to a list of primes only for the first one million integers (i.e. only up to the integer $10^6$), in other words, access to the first 78,498 primes only. In order to obtain digital configurations for higher integers, the expression of the density of primes $1/\ln(x)$ shall be applied, thus enabling us to examine with very high accuracy estimates of the digital proportions within adjacent IPOT intervals of high primes.

Each $(10^N, 10^{(N+1)})$ interval shall be subdivided into nine equally-spaced sub-intervals according to where each digit d leads, namely into the following nine sub-intervals:
$(1*10^N, 2*10^N)$, $(2*10^N, 3*10^N)$, $(3*10^N, 4*10^N)$, … , $(9*10^N, 10*10^N)$.
The prime density for each sub-interval shall be evaluated in its middle, namely at $\{1.5, 2.5, 3.5, 4.5, 5.5, 6.5, 7.5, 8.5, 9.5\}*10^N$. Using the generic expression discussed earlier [**Number of primes within an interval] = [interval's width]*[density**], the number of primes within say $(1*10^N, 2*10^N)$ for digit 1 is obtained as in:
[interval's width]*[density] = $[2*10^N - 1*10^N]*[1/\ln(1.5*10^N)] = 10^N/\ln(1.5*10^N)$.

In general, for any digit d, the number of primes within the sub-interval is:

$[(d+1)*10^N - (d)*10^N]*[1/\ln((d+1/2)*10^N)] =$

$[1*10^N]*[1/\ln((d+1/2)*10^N)] =$

## $10^N / \ln((d + 1/2)*10^N)$

Now we have an almost perfectly precise tool to obtain the digital configurations on any IPOT interval for large primes. Figure 10 depicts digit distributions for a variety of IPOT intervals, omitting the first one from 1 to 10. The first 5 rows in Figure 10, namely from 10 to $10^6$, are the actual digit distributions. From the 6th row on, namely from $10^6$ and up, the distributions are the estimates from $10^N/\ln((d + 1/2)*10^N)$. Hence Figure 10 is comprised of actual and estimated distributions. The high accuracy of the estimates motivated this mixing. For example, for the interval $(10^5, 10^6)$ results are:

Actual:     {12.2%, 11. 6%, 11.4%, 11.1%, 11.0%, 10.8%, 10.8%, 10.6%, 10.5%}
Estimate:   {12.2%, 11.7%,  11.4%, 11.1%, 11.0%, 10.8%, 10.7%, 10.6%, 10.5%}

Namely that with one decimal point precision, and considering primes over $10^6$, digital proportions are practically identical for actual and estimate values. The higher the IPOT interval considered, the closer is this fantastic agreement between actual and estimate, hence all the values for rest of Figure 10 over $10^6$ can be taken as actual ones as well. In other words: everything in the entire table of Figure 10 is actual! Everything is real!

|  |  | 1 | 2 | 3 | 4 | 5 | 6 | 7 | 8 | 9 |
|---|---|---|---|---|---|---|---|---|---|---|
| [ 10 , | 100 ] | 19.0% | 9.5% | 9.5% | 14.3% | 9.5% | 9.5% | 14.3% | 9.5% | 4.8% |
| [ 100 , | 1000 ] | 14.7% | 11.2% | 11.2% | 11.9% | 9.8% | 11.2% | 9.8% | 10.5% | 9.8% |
| [ 1000 , | 10000 ] | 12.7% | 12.0% | 11.3% | 11.2% | 10.7% | 11.0% | 10.1% | 10.4% | 10.6% |
| [ $10^4$ , | $10^5$ ] | 12.4% | 11.8% | 11.5% | 11.1% | 11.0% | 10.5% | 10.8% | 10.5% | 10.5% |
| [ $10^5$ , | $10^6$ ] | 12.2% | 11.6% | 11.4% | 11.1% | 11.0% | 10.8% | 10.8% | 10.6% | 10.5% |
| [ $10^6$ , | $10^7$ ] | 12.0% | 11.6% | 11.3% | 11.1% | 11.0% | 10.9% | 10.8% | 10.7% | 10.6% |
| [ $10^7$ , | $10^8$ ] | 11.9% | 11.5% | 11.3% | 11.1% | 11.0% | 10.9% | 10.8% | 10.7% | 10.7% |
| [ $10^8$ , | $10^9$ ] | 11.8% | 11.5% | 11.3% | 11.1% | 11.0% | 10.9% | 10.9% | 10.8% | 10.7% |
| [ $10^9$ , | $10^{10}$ ] | 11.7% | 11.4% | 11.3% | 11.1% | 11.0% | 11.0% | 10.9% | 10.8% | 10.8% |
| [ $10^{10}$ , | $10^{11}$ ] | 11.7% | 11.4% | 11.2% | 11.1% | 11.0% | 11.0% | 10.9% | 10.8% | 10.8% |
| [ $10^{11}$ , | $10^{12}$ ] | 11.6% | 11.4% | 11.2% | 11.1% | 11.0% | 11.0% | 10.9% | 10.9% | 10.8% |
| [ $10^{12}$ , | $10^{13}$ ] | 11.6% | 11.4% | 11.2% | 11.1% | 11.1% | 11.0% | 10.9% | 10.9% | 10.9% |
| [ $10^{13}$ , | $10^{14}$ ] | 11.5% | 11.3% | 11.2% | 11.1% | 11.1% | 11.0% | 11.0% | 10.9% | 10.9% |
| [ $10^{14}$ , | $10^{15}$ ] | 11.5% | 11.3% | 11.2% | 11.1% | 11.1% | 11.0% | 11.0% | 10.9% | 10.9% |
| [ $10^{15}$ , | $10^{16}$ ] | 11.5% | 11.3% | 11.2% | 11.1% | 11.1% | 11.0% | 11.0% | 10.9% | 10.9% |
| [ $10^{16}$ , | $10^{17}$ ] | 11.5% | 11.3% | 11.2% | 11.1% | 11.1% | 11.0% | 11.0% | 10.9% | 10.9% |
| [ $10^{17}$ , | $10^{18}$ ] | 11.4% | 11.3% | 11.2% | 11.1% | 11.1% | 11.0% | 11.0% | 11.0% | 10.9% |
| [ $10^{18}$ , | $10^{19}$ ] | 11.4% | 11.3% | 11.2% | 11.1% | 11.1% | 11.0% | 11.0% | 11.0% | 10.9% |
| [ $10^{19}$ , | $10^{20}$ ] | 11.4% | 11.3% | 11.2% | 11.1% | 11.1% | 11.0% | 11.0% | 11.0% | 10.9% |
| [ $10^{20}$ , | $10^{21}$ ] | 11.4% | 11.3% | 11.2% | 11.1% | 11.1% | 11.0% | 11.0% | 11.0% | 11.0% |
| [ $10^{21}$ , | $10^{22}$ ] | 11.4% | 11.3% | 11.2% | 11.1% | 11.1% | 11.0% | 11.0% | 11.0% | 11.0% |
| [ $10^{143}$ , | $10^{144}$ ] | 11.2% | 11.1% | 11.1% | 11.1% | 11.1% | 11.1% | 11.1% | 11.1% | 11.1% |
| [ $10^{144}$ , | $10^{145}$ ] | 11.1% | 11.1% | 11.1% | 11.1% | 11.1% | 11.1% | 11.1% | 11.1% | 11.1% |

**Figure 10:** Digital Distributions of Prime Numbers within IPOT Intervals

The table in Figure 10 depicts a mild "digital development pattern" of sorts for the prime numbers, (conceptually) akin to the dramatic Digital Development Pattern observed for all random Benford data where digits start out with an approximate digital equality on the left for low values, continue in the center with roughly the logarithmic configuration, and end up on the far right for high values with extreme digital inequality, even skewer and more severe than that of the Benford inequality, where digit 1 typically usurps leadership by taking over 40% or even over 50% proportion (Kossovsky (2014) chapters 24, 33, 43, 82, 83, and 84.)

# SECTION 3: Density of the Logarithms of Prime Numbers
---

An expression for the 'density' of log of primes shall now be derived applying the Prime Number Theorem. As with $1/\ln(x)$ 'density' for primes which is not a full-fledge density as in mathematical statistics, the density of log of primes to be derived here is also not a full-fledge density in the usual sense, rather it is simply defined as the number of existing log(primes) within 1 unit of log-axis. Surely, if that unity range on the log-axis happened to be bordered exactly by two integers, such as (2, 3) for example, then its mirror-image on the x-axis is an IPOT interval, such as (100, 1000) for the example above.

Let us consider the set of all primes residing within the sequence of consecutive integers from [ x ] to [ x + M - 1] inclusively, where M is a positive whole number.

The estimated number of primes within the integers from [ x ] to [ x + M - 1] is simply the local density of primes times the number of integers in this range (its length).
It should be noted that integer 5 for example occupies one whole unit on the x-axis, namely from 5.0 to 6.0, hence the range of consecutive integers from [ x ] to [ x + M - 1] is from the **left**-most corner of integer x to the **right**-most corner of integer x + M - 1, namely the length of the real x-axis interval (x, x + M) which is x + M – x, or simply M.

Prime Count =  Integer Range * Prime Density
Prime Count = [x + M – x] * (1/ ln(x + M/2))
Prime Count = M/ln(x + M/2)

The approximate midpoint x + M/2 standing roughly between x and x + M is used to evaluate the prime density $1/\ln(t)$.

The number of log(primes) from log(x) to log(x + M) on the log-axis is the same as in the above expression, namely M/ln(x + M/2). Hence Log Density of Primes is:

Log Density of Primes = Count of log(Primes) / Range on the log-axis
Log Density of Primes = M/ln(x + M/2) / [ log(x + M) - log(x) ]
Log Density of Primes = M/ln(x + M/2) /  log((x + M)/(x))
Log Density of Primes = M/ln(x + M/2) /  log(1 + M/x)

$$\text{Log Density of Primes} = \frac{M}{\ln(x + M/2) * \log(1 + M/x)}$$

$$\text{Log Density of Primes} = \frac{x}{\ln\left(x + \frac{M}{2}\right) * \left(\frac{x}{M}\right) * \log\left(1 + \frac{M}{x}\right)}$$

$$\text{Log Density of Primes} = \frac{x}{\ln\left(x + \frac{M}{2}\right) * \log\left(\left(1 + \frac{M}{x}\right)^{\frac{x}{M}}\right)}$$

$$\text{Log Density of Primes} = \frac{x}{\ln\left(x + \frac{M}{2}\right) * \log\left((1 + 1/(\frac{x}{M}))^{\frac{x}{M}}\right)}$$

Now we let M go to zero from above to get the log density in the very immediate vicinity of x. Euler's classic result $\lim_{N \to \infty} \left(1 + \frac{1}{N}\right)^N = e$ enables us to determine that in the limit as M goes to zero - and by implication as x/M goes to infinity - log density is:

$$\text{Log Density of Primes} = \frac{x}{\ln(x) * \log(e)}$$

$$\text{Log Density of Primes} = \frac{x}{\left[\frac{\log(x)}{\log(e)}\right] * \log(e)}$$

$$\text{Log Density of Primes} = \frac{x}{\log(x)}$$

Clearly density of log of primes is monotonically increasing as seen from the expression above [being proportional to x, and inversely proportional only to log(x)], hence the observed near/total digital equality of the sequence of the prime numbers themselves.

Admittedly, for a very small range of integers M, the prime density 1/ln(x) is not so meaningful, as the random overtakes the deterministic and chaos rules over order, much less so for M = 1 or for M being a fraction; for how can one concoct a whole prime out of a fraction of an integer?! A fraction of an integer can never be a whole prime! In any case, one must always keep in mind that the so called 'density' here is simply the <u>fraction</u> of integers that are primes. Yet, letting M approach zero in the limit is an appropriate mathematical procedure here since the focus is temporarily being shifted to the nature of the abstract curves, curves or 'densities' from which the seemingly random occurrences of primes can be deduced and predicted.

Interestingly, even though 1/ln(x) is not a proper density probability function, nonetheless if the Monotonic Transformation Technique in mathematical statistics is used here, the results are identical. The notational convention here is as in John Freund book "Mathematical Statistics", 6th Edition, Section 7.3, Theorem 7.1.

$f(x) = 1/\ln(x)$
$y = u(x) = \log(x)$
$x = w(y) = 10^y$
$g(y) = f[w(y)] * |w'(y)|$
$g(y) = [1/\ln(10^y)] * |(10^y)'|$
$g(y) = [1/\ln(10^y)] * |10^y * \ln(10)|$
$g(y) = [1/\ln(x)] * x * \ln(10)$
$g(y) = [1/[\log(x)/\log(e)]] * x * [\log(10)/\log(e)]$
$g(y) = [\log(e)/\log(x)] * x * [1/\log(e)]$
$g(y) = \mathbf{x/\log(x)}$

Figure 11 depicts the good fit between the theoretical x/log(x) density of log of primes and empirical density of log of primes for the crude and narrow range of 27 consecutive integers, examining 898 primes, from the prime 47, to the prime 7109.

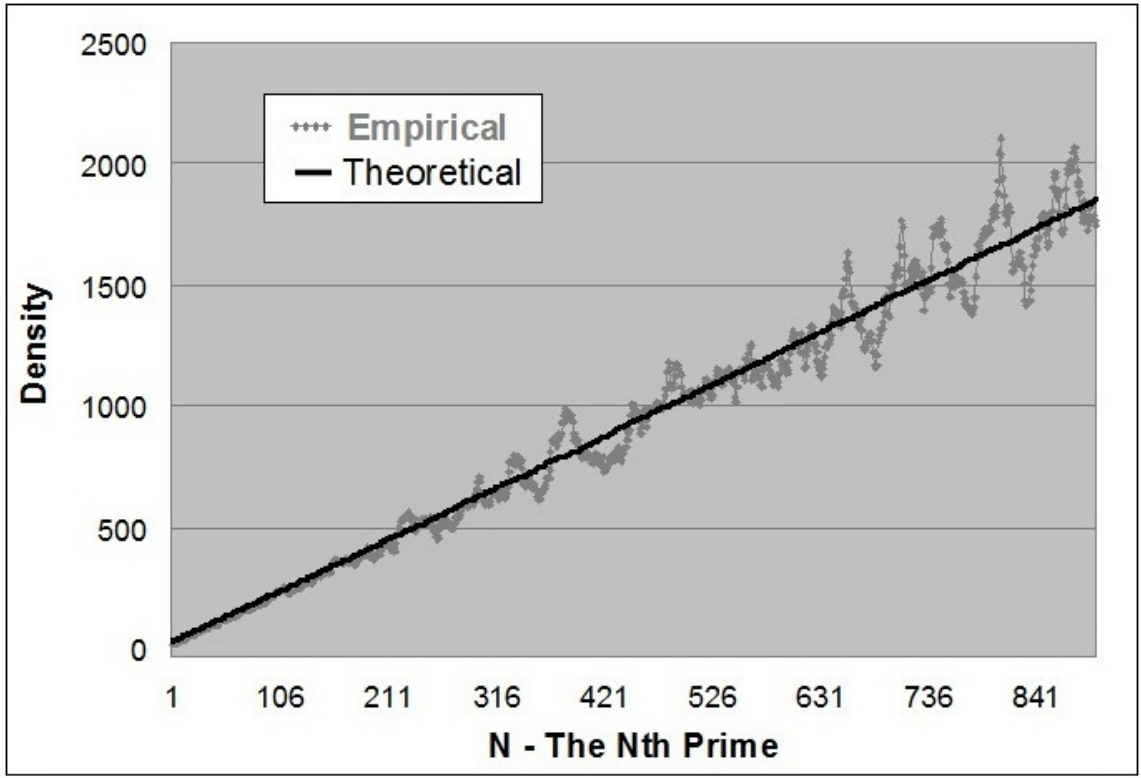

**Figure 11:** Fit of Empirical Log Density and Theoretical x/log(x)

The expression x/log(x) for the density of log of primes is defined as the number of log of primes per 1 unit of log-axis range; hence in essence, it means the number of primes per IPOT interval. Let us evaluate this density expression in the middle of each of the first 6 IPOT intervals from 2 to $10^6$ log-wise, namely at 0.5, 1.5, 2.5, 3.5, 4.5, 5.5. The results show a reasonable fit of theoretical with empirical results as follow:

Theoretical count of primes for the six intervals:  {6, 21, 126, 904, 7027, 57496}
Empirical count of primes for the six intervals:   {4, 21, 143, 1061, 8363, 68906}

For example, for the interval (100, 1000), mid log value is 2.5, corresponding to $10^{2.5}$ or 316.2. Hence log density is evaluated as 316.2/log(316.2), or **126**. This is not too far from the actual number of **143** primes that happened to fall within (100, 1000).

# SECTION 4 :  Dirichlet Density for Subsets of Primes
----------------------------------------------------------------------------------------

Interestingly, even though prime numbers are not Benford, the logarithmic proportions LOG(1 + 1/d) nonetheless pops out here indirectly, as if accidently, not regarding the proportions of numbers with leading digit d, but in another very different sense.

In mathematics, the **Dirichlet density** (also known as the **analytic density**) of a subset of the prime numbers, is a measure of the size of that subset (in comparison to the size of the entire set of all the primes). This measure is easier to use than the natural density.

If A is a subset of the prime numbers, the Dirichlet density of A is the limit (if it exists) of the following ratio as the power s approaches 1 from above:

$$\lim_{s \to 1+} \frac{\sum_{p \in A} \frac{1}{p^s}}{\sum_{p} \frac{1}{p^s}}$$

In other words:

$$\lim_{s \to 1+} \frac{\sum_{\text{only the primes in subset } A} \frac{1}{p^s}}{\sum_{\text{all the primes}} \frac{1}{p^s}}$$

Clearly, the denominator contains the numerator; the numerator is less than or equal the denominator; so that this limit-ratio - if it exists - is less than or equal to 1, and at a minimum it is zero, but it is never negative. Hence Dirichlet density is in [0, 1].

For example, the Mersenne subset of the prime numbers are those primes that are one less than an integral power of two, namely $2^N - 1$, N being an integer. The first seven **Mersenne primes** are $\{2^2 - 1, 2^3 - 1, 2^5 - 1, 2^7 - 1, 2^{13} - 1, 2^{17} - 1, 2^{19} - 1\}$, namely: {3, 7, 31, 127, 8191, 131071, 524287}. They appear to be quite rare, thinly spread among the entire set of all primes. The Dirichlet density for Mersenne primes is then:

$$\lim_{s \to 1+} \frac{\frac{1}{3^s} + \frac{1}{7^s} + \frac{1}{31^s} + \frac{1}{127^s} + \frac{1}{8191^s} + \cdots \text{ to } \infty}{\frac{1}{2^s} + \frac{1}{3^s} + \frac{1}{5^s} + \frac{1}{7^s} + \frac{1}{11^s} + \frac{1}{13^s} + \frac{1}{17^s} + \frac{1}{19^s} + \cdots \text{ to } \infty}$$

The value of this limit-ratio is zero, indicating conceptually that there are very few Mersenne primes within the entire set of the primes and that they are quite rare indeed. Even though we are dealing with concepts such as hierarchy of infinities and a limit process, a very informal or simple-minded assessment of the rarity of Mersenne primes within all the primes can be given by the extremely low percentage of their occurrences up to the finite integer 524287, namely (7)/(524287), or merely 0.0013%!

## SECTION 5:  Modified Dirichlet Density
-----------------------------------------------------------------

As another example, assuming one could obtain the limit as s goes to 1 from above by direct substitution of 1 for s, then for the consideration of the subset of all the primes with 3 as the first leading digit, the Dirichlet density is calculated simply as:

$$\frac{\frac{1}{3} + \frac{1}{31} + \frac{1}{37} + \frac{1}{307} + \frac{1}{311} + \frac{1}{313} + \cdots \; to \; infinity}{\frac{1}{2} + \frac{1}{3} + \frac{1}{5} + \frac{1}{7} + \frac{1}{11} + \frac{1}{13} + \frac{1}{17} + \frac{1}{19} + \cdots \; to \; infinity}$$

The informal and simplistic expectation here is of a non-zero ratio given the substantial 11.1% proportion approximately that primes between integral powers of ten occur having 3 as their first digit [i.e. eventual digital equality within IPOT intervals.]

But direct substitution of 1 for s in the above estimation of the limit is a suspect, since it has not been shown mathematically that plugging in 1 for s is actually the limit itself. Let us demonstrate the challenges and issues involved in attempting to prove this:

The division law of limits:

If the limits of $\lim_{x \to c} f(x)$ and $\lim_{x \to c} g(x)$ both exist, and limit of $\lim_{x \to c} g(x) \neq 0$, then:

$$\lim_{x \to c} \frac{f(x)}{g(x)} = \frac{\lim_{x \to c} f(x)}{\lim_{x \to c} g(x)}$$

The addition law of limits:

If the limits of $\lim_{x \to c} f(x)$ and $\lim_{x \to c} g(x)$ both exist, then:

$$\lim_{x \to c} [f(x) + g(x)] = [\lim_{x \to c} f(x)] + [\lim_{x \to c} g(x)]$$

The Exponential Function $q(x) = B^x$ is continuous at a minimum for all $x \geq 0$ and for all $B \geq 1$ (a condition that all prime numbers satisfy within the repeated $p^s$ expressions of the Dirichlet density, since $p_j \geq 2$). Hence the limit of $B^x$ as x approaches 1 (from above or from below and in general) certainly exists, and it can be obtained by direct substitution. $B^x$ is also always a positive non-zero quantity for all those B and x values. From the division law of limits it follows that for the function $h(s) = \dfrac{1}{p^s}$ where $p \geq 1$ and $s \geq 0$,

$$\lim_{s \to 1+} \frac{1}{p^s} = \frac{1}{p^1}$$

; namely that direct substitution is allowed; and that this limit certainly exists. Consequently, by the addition law of limits, each part of the Dirichlet density separately, namely the numerator as well as the denominator, can be evaluated stepwise, via infinitely many direct substitutions of 1 for s. Yet, unfortunately, the limit of the numerator as well as the limit of the denominator do not exist, since each one diverges (adding infinitely many 1/p terms which never approaches zero close enough yields an ever increasing sum, namely it's infinite). Surely the denominator is decisively a positive non-zero quantity since all prime numbers are positive and there are only addition terms within the Dirichlet expression, not involving any subtractions. Yet, the division law of limits cannot be utilized here for the limit of main division in the Dirichlet density since it's in an indeterminate form where both the numerator and the denominator diverge. Unfortunately, L'Hôpital's rule does not help here either.

In any case, the above simplification – from a limiting process to direct substitution of 1 for s – whether it is justified mathematically in the infinite Dirichlet case or not, shall be utilized in the definition (*hence correctly*) of a modified finite form of Dirichlet in order to gain better theoretical understanding of the forces at play here leading to the logarithmic proportions.

In addition, the quest to perform an approximate computerized test for the Dirichlet density is frustrated since the expression involves an infinite number of primes in the denominator as well as an infinite number of primes with d as the leading digit in the numerator. The suggested next step then in modifying Dirichlet density is to perform partial Dirichlet densities with s = 1 in a piecemeal and finite manner, one integral-power-of-ten-interval at a time, considering a ratio pertaining to two subsets, that of the denominator being the set of all the primes in $(10^{INTEGER}, 10^{INTEGER + 1})$, and that of the numerator being the smaller subset of all the primes with d as the leading digit in $(10^{INTEGER}, 10^{INTEGER + 1})$. Such an approach would enable us to examine how Dirichlet density works locally at each IPOT interval level. This ratio shall be termed 'The Modified Dirichlet Density'.

MDD(N, d) ≡ Modified Dirichlet Density(N, d) ≡

$$\frac{\sum_{\text{primes with leading digit } d \text{ in } (10^N, 10^{N+1})} \frac{1}{p}}{\sum_{\text{all the primes in } (10^N, 10^{N+1})} \frac{1}{p}}$$

It is hoped that a steady and clear pattern would emerge in these proportions across all IPOT intervals, a pattern which could serve as an indication of the overall proportions for all the prime numbers between 2 and infinity, namely regarding the original all inclusive definition of Dirichlet Density.

Having access to the actual set of all the primes from 2 to 1,000,000 enables the author to perform such calculations for five IPOT intervals. The initial interval of (2, 10) is omitted as being an outlier of sorts, not having enough primes to show any stable pattern.

Empirical results for Modified Dirichlet Density for N values {1, 2, 3, 4, 5} are as follows (one detailed example of the calculations of **44.6%** is given in next page):

Primes 10 to 100            {44.6, 12.4, 9.5, 11.0, 5.7, 5.0, 6.5, 3.7, 1.6}   SSD = 265.6
Primes 100 to 1000          {37.0, 16.1, 11.6, 9.6, 6.4, 6.3, 4.7, 4.5, 3.7}   SSD = 55.1
Primes 1000 to 10,000       {33.1, 18.1, 12.1, 9.3, 7.3, 6.3, 5.0, 4.5, 4.1}   SSD = 11.2
Primes 10,000 to 100,000    {32.3, 18.0, 12.4, 9.3, 7.6, 6.1, 5.4, 4.6, 4.2}   SSD = 6.3
Primes 100,000 to 1,000,000 {32.0, 17.9, 12.4, 9.4, 7.6, 6.3, 5.4, 4.7, 4.2}   SSD = 4.6
**Benford's Law first digits:   {30.1, 17.6, 12.5, 9.7, 7.9, 6.7, 5.8, 5.1, 4.6}   SSD = 0**

[Note: The vector of nine elements refers to the values of Modified Dirichlet Density starting with digit 1 on the left, and ending with digit 9 on the right. These values are actually ratios, but thought of as percent, while the '%' sign is omitted for brevity.]

**Eureka!** An apparent pattern converging to Benford is found <u>locally</u> on each finite interval between adjacent IPOT.

It is reasonable to conjecture that this pattern should persist for higher primes, consistently and gradually approaching the logarithmic proportions in the limit as primes go to infinity.

Remarkably, this peculiar calculation for each of the 9 digits leads directly to the Benford proportions of LOG(1 + 1/d). It must be emphasized that we are not calculating proportions of primes with d as the leading digit, but rather the sum of the reciprocals of primes with digit d leading, divided by the sum of the reciprocals of all those primes regardless of which d is leading, and all within one IPOT interval. Clearly, this is <u>not</u> just another manifestation of Benford's Law, but rather just a curious coincidence perhaps.

It would be a fallacy to state that: "**The prime numbers are Benford Dirichlet-wise**". They are not! It is only this particular expression concocted out of the primes that somehow happened to be LOG(1 + 1/d).

One detailed example of the empirical calculations above:

The entire set of all the primes between 10 and 100 is:
{11, 13, 17, 19, 23, 29, 31, 37, 41, 43, 47, 53, 59, 61, 67, 71, 73, 79, 83, 89, 97}.

Our modified piecemeal Dirichlet density for digit 1 is calculated as follows:

$$\text{MDD}(1, 1) = \frac{\Sigma \; primes \; in \; (10,100) \; with \; 1 \; as \; the \; leading \; digit \; \frac{1}{p}}{\Sigma \; all \; the \; primes \; in \; (10,100) \frac{1}{p}}$$

$$\frac{\frac{1}{11} + \frac{1}{13} + \frac{1}{17} + \frac{1}{19}}{\frac{1}{11} + \frac{1}{13} + \frac{1}{17} + \frac{1}{19} + \frac{1}{23} + \; ... \; + \frac{1}{83} + \frac{1}{89} + \frac{1}{97}}$$

$$\frac{0.0909 + 0.0769 + 0.0588 + 0.0526}{0.0909 + 0.0769 + 0.0588 + 0.0526 + 0.0435 + \; ... \; + 0.0120 + 0.0112 + 0.0103}$$

$$\frac{0.2793}{0.6266} = 0.446 = 44.6\%$$

# SECTION 6: Explaining Modified Dirichlet Connection to Benford
--------------------------------------------------------------------------------------------------------

In order to take the mystery out of this Modified Dirichlet surprising result, let us examine the actual quantity being calculated here, part of which is simply the sum of reciprocals of primes, first digit by first digit. Little reflection is needed to realize that this should yield a decisive advantage for low digits. This follows from the fact that within any particular interval bounded by adjacent IPOT points, a number with a low leading digit implies low quantity, and a number with a high leading digit implies high quantity, namely that leading digits and actual quantities positively correlate (Kossovsky (2014), chapter 143 titled "Digits Serving as Quantities in Benford's Law").

Hence on the IPOT interval (10, 100): the reciprocals of the primes **11** or **23** (with low first digits 1 and 2) yield decisively higher values than the reciprocals of the primes **83** or **97** (with high first digits 8 and 9).

And on the IPOT interval (100, 1000): the reciprocals of the primes **101** or **211** (with low first digits 1 and 2) yield decisively higher values than the reciprocals of the primes **887** or **997** (with high first digits 8 and 9).

Primes within IPOT intervals start out with slightly skewed digital proportions in favor of low digits, but are not too far from equality. In the limit going to infinity, primes on such intervals end up with digital equality, namely that the local histogram or density of primes at infinity is uniform with a flat curve. Hence with no defense given to high digits say in the form of occurring more often than low digits, merely summing reciprocals along digital lines implies a clear advantage for low digits. The only question is by how much exactly do low digits gain over high digits in this reciprocal game?

Between 100 and 1000, the reciprocal of the prime 151 for example is 1/151 or 0.0066, while the reciprocal of the prime 953 for example is 1/953 or 0.0010. Hence the reciprocal of 151 is **6.31** times as large as the reciprocal of 953! Adding slightly to this dichotomy is the fact that early on (before we hit infinity) lower digits occur a bit more frequently (14.7% for digit 1 on the interval [100, 1000]) than high digits (9.8% for digit 9 on the interval [100, 1000]), hence in this example, low digit 1 occurs approximately (14.7)/(9.8) or **1.50** times more frequently than high digit 9. The confluence of these two factors [higher reciprocals and higher frequencies] equals roughly (6.31)*(1.50), or 9.46, being the factor by which reciprocals of digit 1 have an advantage over reciprocals of digit 9, and which is quite close to the actual digit 1 to digit 9 ratio (37.0)/(3.7) or 10.0 of the Modified Dirichlet density empirically calculated above. The classic Benfordian advantage is LOG(1 + 1/1)/LOG(1 + 1/9), or (30.1)/(4.6), namely that digit 1 occurs **6.58** times more frequently than digit 9.

The above discussion justifies the expectation that the pattern of the empirical results obtained earlier for the Modified Dirichlet Density for N values {1, 2, 3, 4, 5} should converge to Benford in the limit as primes go to infinity. This is so since ultimately at infinity all nine digits attain an equal proportion of 1/9 or 11.1% within any single IPOT interval, and thus Modified Dirichlet Density approaches the classic Benfordian digit-1-over-digit-9 advantage of the 6.58 factor (as shall be seen further here in details). In other words, that the driving force behind the gradual digital development of the Modified Dirichlet, from the beginning with results that are slightly skewer than the logarithmic, giving digit 1 a bit more than 30.1%, and into the eventual Benford proportion, is all about digits themselves developing from slight skewness in the beginning favoring low digits, and into equality in the limit. Once digits 'settle into' their equal proportions, Modified Dirichlet Density attains its ultimate logarithmic configuration.

Let us now demonstrate in more details the mechanism at play here, and then prove that Modified Dirichlet Density should approach $LOG(1 + 1/d)$ in the limit as primes go to infinity. We first start with the particular example of all the primes on (10, 100), in order to demonstrate visually how Modified Dirichlet relates indirectly to Benford, and then prove the result in general for any IPOT interval as primes go to infinity.

Figure 12 depicts the reciprocals of all the primes on (10, 100), namely of the 21 primes {11, 13, 17, 19, 23, 29, 31, 37, 41, 43, 47, 53, 59, 61, 67, 71, 73, 79, 83, 89, 97}.

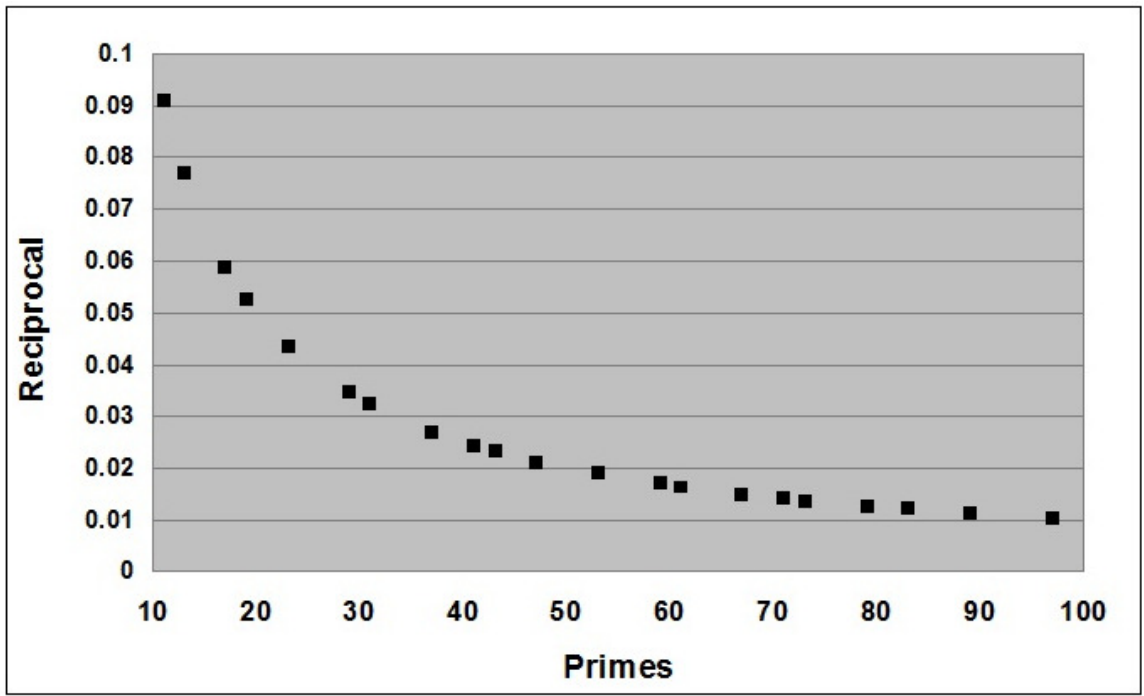

**Figure 12:** The Reciprocals of the 21 Primes on (10, 100)
`

$$\text{MDD}(1, 1) = \frac{\sum \textit{primes with 1 as leading digit in } (10,\ 100)\ \frac{1}{p}}{\sum \textit{all the primes in } (10,\ 100)\ \frac{1}{p}}$$

Let us multiply the numerator and the denominator by the average distance between primes for this particular interval between 10 and 100 IPOT points. Calculating the average distance here yields the value 4.3. We get:

$$\text{MDD}(1, 1) = \frac{Avg\ Dist * \sum \textit{primes in } (10,\ 100)\ 1\ leads\ \frac{1}{p}}{Avg\ Dist * \sum \textit{all the primes in } (10,\ 100)\ \frac{1}{p}}$$

$$\text{MDD}(1, 1) \approx \frac{Area\ of\ Reciprocals\ of\ Digit\ 1}{Entire\ Area\ of\ Reciprocals\ of\ all\ Digits}$$

$$\text{MDD}(1, 1) \approx \frac{Shaded\ Area\ in\ Figure\ 13}{Entire\ Area\ in\ Figure\ 14}$$

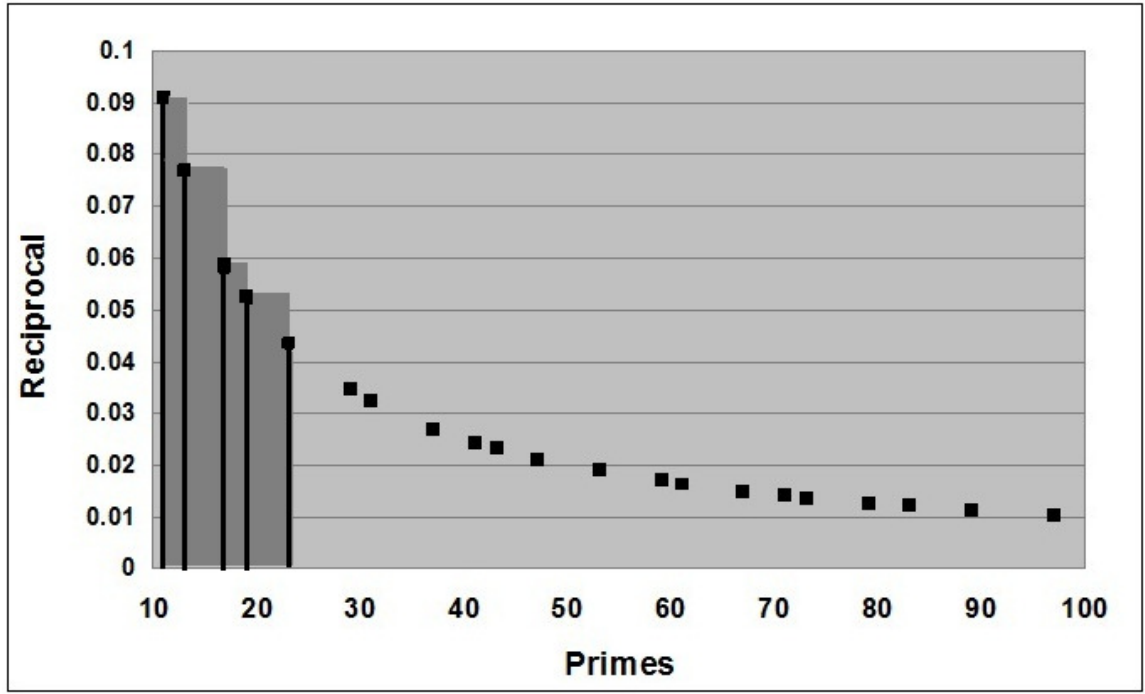

**Figure 13:** Area for the Reciprocals of Primes with digit 1 Leading

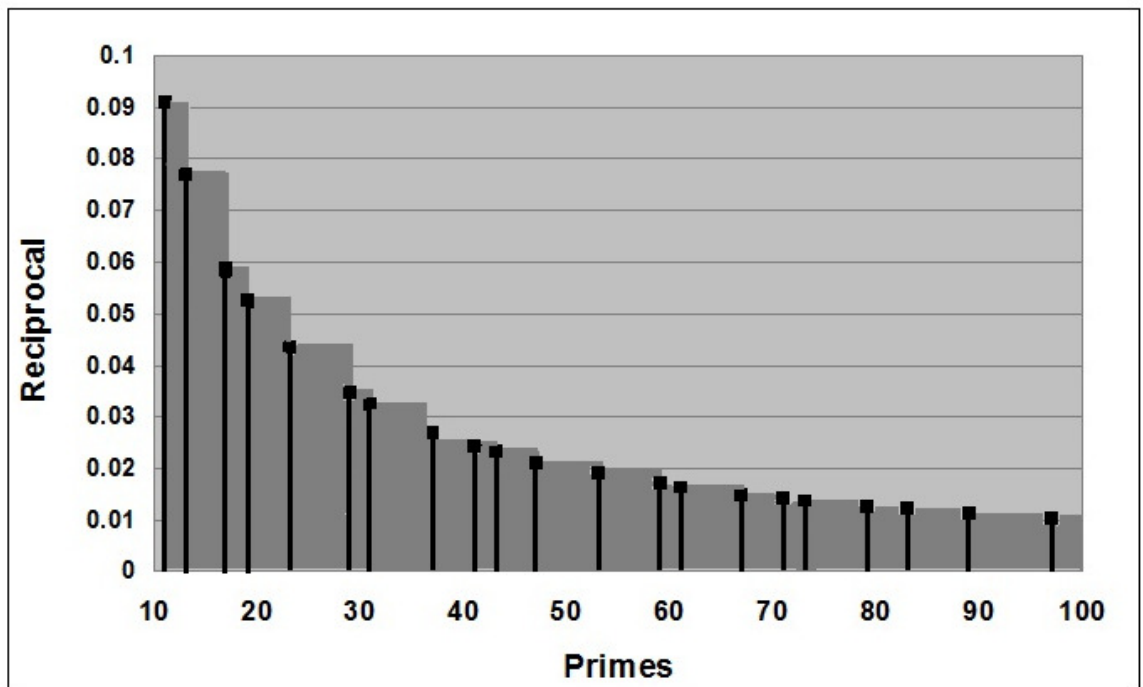

**Figure 14:** The Entire Area for Reciprocals of All the Primes in Interval

Without loss of generality, the drawing in Figures 13 and 14 utilize the actual distances between each pair of adjacent primes, instead of using the average distance 4.3. Does it matter much? Certainly there are offsetting and compensating errors here; some rectangles are too short [in relation to the average distance]; others are too long; but overall the errors should cancel each other out; especially for much higher IPOT intervals where the primes are much more numerous, and are more smoothly distributed, as shall be shown shortly.

In order to visually drive the point that low digits always win the reciprocal game within IPOT intervals, Figure 15 is added showing the area of reciprocals for high digit 7, which is much more modest in comparison to area of reciprocals for low digit 1 of Figure 13.

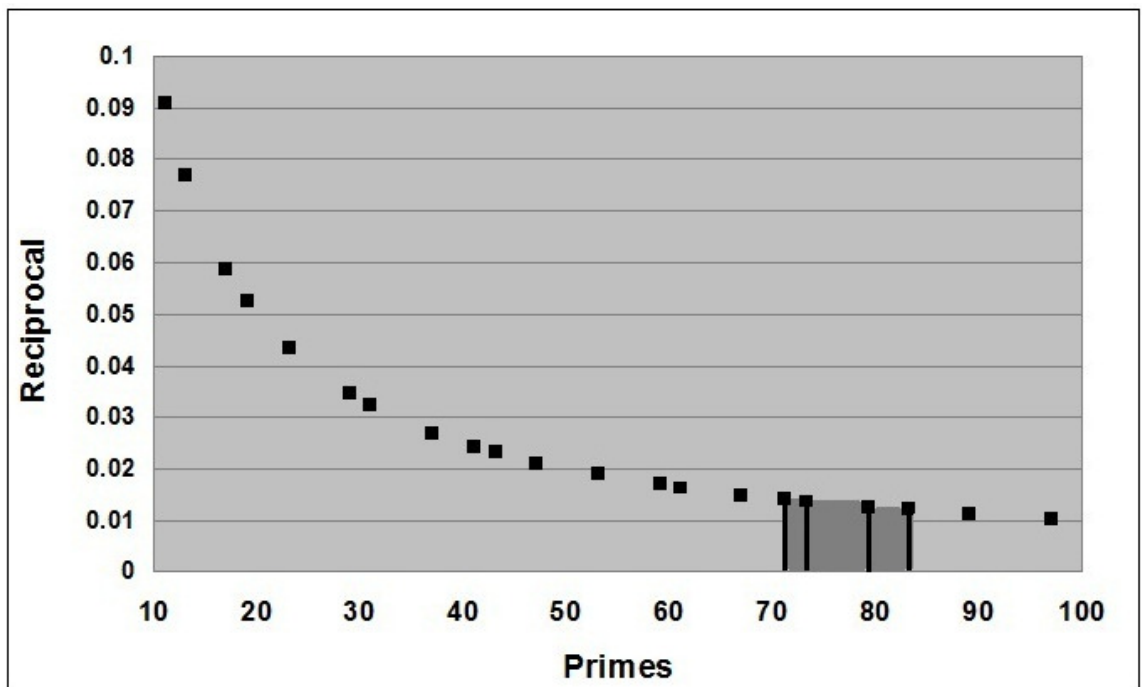

**Figure 15:** Area for the Reciprocals of Primes with digit 7 Leading

Surely, a differential in concentration of primes, not being approximately uniformly spread out at the macro level, would invalidate the above model which converts Modified Dirichlet Density into ratio of two areas (via the singular substitution of the average distance for all actual distances between primes). In reality, at the macro level, distances between primes do increase as higher integers are considered, and this in essence is the Prime Number Theorem; that primes become sparse and rarer as larger integers are considered; but within any given IPOT interval this increase in distance is quite marginal. Moreover, in the limit as primes go to infinity, distances between primes are approximately equal within any particular IPOT interval, therefore the use of the average distance between primes throughout as a singular quantity, for the left part of the interval around digit 1, for the center around digit 5, and for the right part of the interval around digit 9, is certainly justified.

Here are (five empirical and one estimated) results regarding the typical distances between primes (in units of an integer of course).

Two results are given for each IPOT interval:

(1) <u>at the beginning</u> of the interval where digit 1 leads;
(2) <u>at the end</u> of the interval where digit 9 leads:

Average distance between primes on (10, 20) is **2.6**
Average distance between primes on (90, 99) is **6.0**

Average distance between primes on (100, 200) is **4.9**
Average distance between primes on (900, 999) is **7.9**

Average distance between primes on (1,000, 2,000) is **7.4**
Average distance between primes on (9,000, 9,999) is **8.7**

Average distance between primes on (10,000, 20,000) is **9.7**
Average distance between primes on (90,000, 99,999) is **11.4**

Average distance between primes on (100,000, 200,000) is **11.9**
Average distance between primes on (900,000, 999,999) is **13.8**

Average distance between primes farther on $(1*10^{37}, 2*10^{37})$ is **85.6**
Average distance between primes farther on $(9*10^{37}, 10*10^{37})$ is **87.4**

The last pair is the highly accurate estimates utilizing the expression derived earlier $10^N/\ln((d + 1/2)*10^N)$ for the number of primes with digit d leading within IPOT interval, namely: Average Distance = $10^N/[10^N/\ln((d + 1/2)*10^N)] = \ln((d + 1/2)*10^N)$. The dramatic/gentle dichotomy above within each early/late interval is simply the mirror-image of the fact that low digits are significantly/slightly more frequent than high digits in the beginning before we hit infinity.

Students of Benford's Law are immediately struck by the strong conceptual similarity between the expression of the Modified Dirichlet Density and the most perfectly logarithmic distribution of them all, namely **PDF(x) = k/x** defined over the range between two adjacent IPOT points. Isn't reciprocity the most essential quality of the k/x distribution?! Could the well-known results about k/x distribution explain MDD? (Kossovsky (2014), chapter 60 titled "The Case of k/x Distribution", and chapter 62 titled "Uniqueness of k/x Distribution".)

In order to complete the proof, it is necessary to formalize the comparison between: (1) discrete sums of reciprocals of primes with digit d leading, and (2) continuous area under k/x curve over the interval $((d)10^N, (d+1)10^N)$ with d running from 1 to 9 signifying first leading digits. Let us now demonstrate this comparison more rigorously.

The occurrences of primes on (10, 100) in Figures 12, 13, 14, and 15, are very particular, as if each prime has its own unique personality; as if distances between them signify something important. As higher IPOT intervals are considered, occurrences of primes appear more smooth, a bit more predictable, and primes lose their identities, becoming more ordinary and repetitive. Figure 16 depicts the reciprocals of all the primes on (100, 1000); visually confirming this generic observation. Figure 16 is visually also strongly suggestive of the k/x density curve, much more so than Figure 12 is. Other prime-reciprocal curves constructed for much larger IPOT intervals should resemble k/x density even better.

Figure 17 depicts the perfectly Benford distribution of k/x on (100, 1000).
In order to equate entire area to one, we set $\int k/x \, dx = 1$ with [100, 1000] as the limits of integration, hence $k*[\ln(1000) - \ln(100)] = 1$, so that $k*\ln(1000/100) = 1$, and finally $k = 1/\ln 10 = \log(e) \approx 0.434294482$. Clearly $\int \log(e)/x \, dx$ with $[(d), (d+1)]$ as the limits of integration yields $\log(e)*[\ln(d+1) - \ln(d)] = \log(1 + 1/d)$, namely Benford's Law.

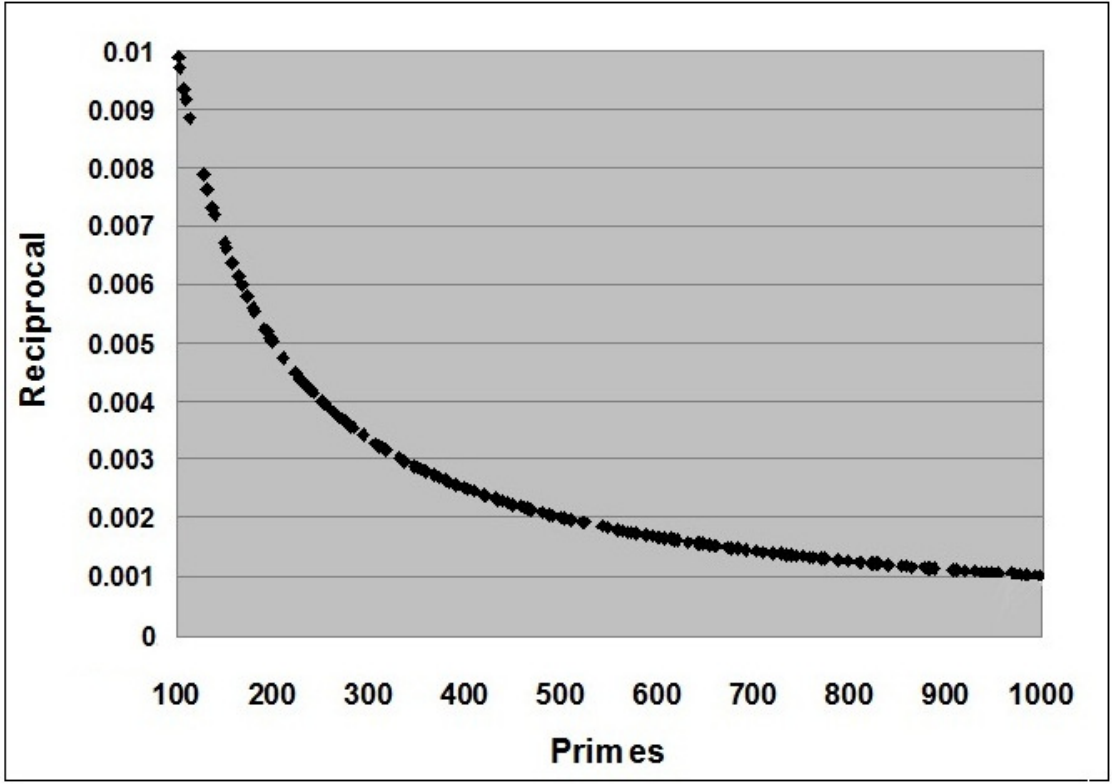

**Figure 16:** Reciprocals of Primes on (100, 1000) Strongly Resemble k/x

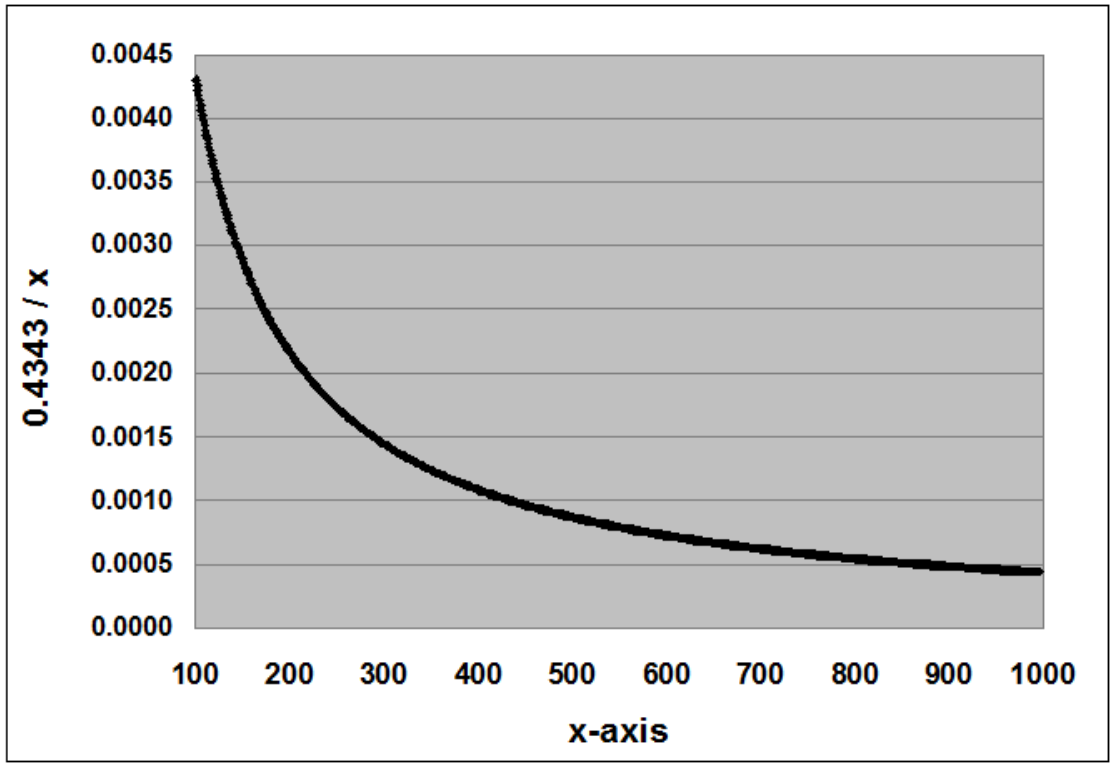

**Figure 17:** The Perfectly Benford Distribution [0.4343]*1/x on (100, 1000)

An additional factor which drives each consecutive prime-reciprocal curve between IPOT points closer towards k/x distribution is that within each new cycle between IPOT there are more primes; they keep growing!

There are      4  primes between 1 and 10
There are     21  primes between 10 and 100
There are    143  primes between 100 and 1,000
There are   1061  primes between 1,000 and 10,000
There are   8363  primes between 10,000 and 100,000
There are  68906  primes between 100,000 and 1,000,000

The expression $10^N/\ln((d + 1/2)*10^N)$ for the 9 digital sections within each IPOT interval obtained earlier from the density of primes $1/\ln(x)$ can be applied in estimating with extremely high accuracy the number of primes within each IPOT interval. An excellent close fit of theoretical with empirical/actual results is obtained as follow:

Empirical count of primes for the six intervals:   {4, 21, 143, 1061, 8363, 68906}
Theoretical count of primes for the six intervals: { 8, 24, 147, 1068, 8379,  68970}

The empirical growth factors are {21/4, 143/21, 1061/143, 8363/1061, 68906/8363}, namely {5.3, 6.8, 7.4, 7.9, 8.2}. Computer calculations of the estimates of the next prime growth factors between IPOT intervals beyond 1,000,000 applying the expression $10^N/\ln((d + 1/2)*10^N)$ for the 9 digital sections within each IPOT interval yields: {8.5, 8.7, 8.8, 9.0, 9.1, 9.1, 9.2, 9.3, 9.3, 9.4, 9.4, 9.4, 9.5, 9.5, 9.5, 9.5, 9.6, etc.}.
In the limit, prime growth factor between consecutive IPOT intervals approaches **10-fold**. This is easily explained under the assumption that density $1/\ln(x)$ for primes barely changes for very large x, attainting approximately the same level on $(10^N, 10^{N+1})$ as on $(10^{N+1}, 10^{N+2})$, hence since range of the latter is 10 times as big as the former, the expression [Prime Count] = [Integer Range]*[Prime Density] implies that there are 10 times as many primes in the latter as compared with the former.

Therefore, as much larger primes are considered, each IPOT interval contains huge number of primes, and overall spread is highly smooth, strongly resembling k/x curve. For that reason the model converting the Modified Dirichlet Density into the ratio of two areas is certainly justified for higher primes on the way to infinity when particular distances between primes can be replaced with confidence by the average distance. Let us complete the proof for any MDD(N, d) with sufficiently high N:

$$\text{MDD}(N, d) = \frac{\sum_{(d)10^N}^{(d+1)10^N} \frac{1}{p}}{\sum_{10^N}^{10^{(N+1)}} \frac{1}{p}}$$

Where only primes are considered within the indices of summation, omitting the composites. Multiplying the numerator and the denominator by the average distance between primes within $(10^N, 10^{(N+1)})$ we get:

$$\text{MDD}(N, d) = \frac{(\text{Avg Dist}) * \sum_{(d)10^N}^{(d+1)10^N} \frac{1}{p}}{(\text{Avg Dist}) * \sum_{10^N}^{10^{(N+1)}} \frac{1}{p}}$$

Interpreting the sums as areas under $1/x$ curve, we write:

$$\text{MDD}(N, d) = \frac{\int_{(d)10^N}^{(d+1)10^N} \frac{1}{x} dx}{\int_{10^N}^{10^{(N+1)}} \frac{1}{x} dx}$$

Evaluating the definite integrals, we get:

$$\text{MDD}(N, d) = \frac{\ln\left((d+1)10^N\right) - \ln\left((d)10^N\right)}{\ln\left(10^{(N+1)}\right) - \ln(10^N)}$$

$$\text{MDD}(N, d) = \frac{\ln(d+1) + \ln(10^N) - \ln(d) - \ln(10^N)}{(N+1)*\ln(10) - N*\ln(10)}$$

$$\text{MDD}(N, d) = \frac{\ln(d+1) - \ln(d)}{N*\ln(10) + \ln(10) - N*\ln(10)}$$

$$\text{MDD}(N, d) = \frac{\ln((d+1)/d)}{\ln(10)} = \frac{\ln(1+1/d)}{\ln(10)} = \frac{\log(1+1/d)}{\log(10)}$$

$$\text{MDD}(N, d) = \frac{\log(1+1/d)}{1} = \text{The Benford Proportions}$$

# SECTION 7:  Modified Dirichlet Leads Original  Dirichlet to Benford
---

Could the above result about the Modified Dirichlet Density tell us something about the original all inclusive Dirichlet Density which considers all the primes in one fell swoop? The answer is in the affirmative; in fact the above result implies logarithmic proportions also for the original Dirichlet Density. The challenging issue here is the profound dichotomy between the two densities; namely that while the modified one had already used direct substitution of 1 for s and now comprises of reciprocals; the original one is still a limiting process as variable s approaches 1 from above.

It was shown earlier that from $10^{145}$ digital proportions of primes are equal considering one decimal precision (see Figure 10). Hence beyond $10^{145}$ one can confidently assume that the Modified Dirichlet Density is nearly Benford. Define k >145 as a sufficiently large power of ten such that from $10^k$ onward MDD(k, d) is Benford and digital equality prevails for a fixed m-decimal-precision criteria.

Let LMDD(k, d), LMDD(k + 1, d), LMDD(k + 2, d), LMDD(k + 3, d), etc. expressed as **LMDD(k + j, d)** be the Modified Dirichlet Density on the intervals $(10^{k+j}, 10^{k+j+1})$ as a limit process with s approaching 1 from above, and with j running as 0, 1, 2, 3, and so forth, namely MDD as the limits of the ratios, where each 1/p term is written as $1/p^s$, and considering all IPOT intervals beyond $10^k$.

Incorporating a limit process in the above definition yields the same result for all MDD ratios as shown earlier via the applications of the addition and division laws of limits, enabling us to evaluate the limits of the individual terms:

$$\lim_{s \to 1+} \frac{1}{p^s} = \frac{1}{p^1}$$

The limit of each LMDD exists since it incorporates only finite terms, namely that the numerator and denominator both converge. Surely LMDD = MDD.

Define **N(j)** as the numerator and **D(j)** as the denominator for each LMDD(k + j, d) involving variable s, namely as:

$$LMMD(k+j, d) = \lim_{s \to 1+} \frac{N(j)}{D(j)}$$

for some fixed d value, although d is not shown for brevity.

The sequence of LMDD(k + j, d) with j = 0, 1, 2, 3, and so on, is then the limit as s goes to 1 from above of:

$$\frac{N(0)}{D(0)}, \frac{N(1)}{D(1)}, \frac{N(2)}{D(2)}, \frac{N(3)}{D(3)}, \ldots, \text{LOG}(1 + 1/d), \text{LOG}(1 +1/d), \ldots$$

The notation B is used for the Benford proportions of LOG(1 +1/d) for the fixed d value in the definition of N(j) and D(j). Hence the following infinitely many relationships hold in the limit as s approaches 1 from above:

$$B = \frac{N(0)}{D(0)} \quad B = \frac{N(1)}{D(1)} \quad B = \frac{N(2)}{D(2)} \quad B = \frac{N(3)}{D(3)} \quad \text{etc. to infinity}$$

Hence:

D(0)*B = N(0)
D(1)*B = N(1)
D(2)*B = N(2)
D(3)*B = N(3)

and so forth to infinity.

The original all inclusive (infinite) Dirichlet Density in essence is built from the various (finite) Limit Modified Dirichlet expressions. For example, for digit 3 it is defined as:

$$\lim_{s \to 1+} \frac{\frac{1}{3^s} + \frac{1}{31^s} + \cdots \frac{1}{307^s} + \cdots \frac{1}{3001^s} + \cdots \frac{1}{30011^s} + \cdots \text{to } \infty}{\frac{1}{2^s} + \frac{1}{3^s} + \frac{1}{5^s} + \frac{1}{7^s} + \frac{1}{11^s} + \frac{1}{13^s} + \frac{1}{17^s} + \frac{1}{19^s} + \frac{1}{23^s} + \cdots \text{to } \infty}$$

And one can find within the above expression all the Limit Modified Dirichlet ones for leading digit 3, without overlapping, without repeating any term. In other words, the above expression is comprised solely and exactly of the various Limit Modified Dirichlet ones for leading digit 3. Ignoring the contributions from earlier intervals below $10^k$ which are negligible as primes go to infinity; as s approaches 1 from above for the original Dirichlet Density, it can be written analytically (finite) term by (finite) term as:

$$\lim_{s \to 1+} \frac{N(0) + N(1) + N(2) + N(3) + \cdots \text{ to infinity}}{D(0) + D(1) + D(2) + D(3) + \cdots \text{ to infinity}}$$

Yet, it is still a limit, not a ratio, since both the numerator and denominator diverge, implying that the limit is in an indeterminate form. Substituting B*D(j) for all the N(j) terms (*true though only when s goes to 1 – but s does indeed go to 1 here!*) we get:

$$\lim_{s \to 1+} \frac{B * D(0) + B * D(1) + B * D(2) + B * D(3) + \cdots \text{ to infinity}}{D(0) + D(1) + D(2) + D(3) + \cdots \text{ to infinity}}$$

$$\lim_{s \to 1+} \frac{B * [D(0) + D(1) + D(2) + D(3) + \cdots \text{ to infinity}]}{D(0) + D(1) + D(2) + D(3) + \cdots \text{ to infinity}}$$

$$\lim_{s \to 1+} \frac{B}{1} = B$$

Hence the limit of the original Dirichlet density is LOG(1 +1/d) and this completes the proof.

# SECTION 8: Conceptual Conclusions
----------------------------------------------------

Prime numbers are decisively non-Benford as observed in the normal density of the proportion of primes with digit d leading, in spite of the Dirichlet results obtained in this article. Primes are innately uniformly distributed with a flat density in the limit having digital equality at infinity. Mr. Johann Dirichlet artificially, arbitrarily, and forcefully, dresses the primes with that reciprocal coat, making them appear as if they are Benford. The primes suffer in silence, pretending to go along, but deep down they continue to believe in uniform, equitable, and fair digital distribution, hoping that one day they would be able to tear down that Dirichlet's unnatural dress and act naturally. In the same vein, Dirichlet might as well employ the same reciprocal coat and dress up various other sequences that have nothing to do with the logarithmic distribution, and then proclaim them all to be Benford!?

Any argument in support of the false statement that "The prime numbers are Benford Dirichlet-wise" pointing to the facts that MDD and DD are LOG(1 + 1/d), should also be countered by pointing to the sequence $X_{N+1} = X_N + 1$ beginning with $X_1 = 1$ (namely the sequence of all the integers Z), and constructing a similar limit-ratio (involving not only the primes but also the composites) such as:

$$\lim_{s \to 1+} \frac{\Sigma \text{ integers with leading digit d } \frac{1}{z^s}}{\Sigma \text{ all the integers } \frac{1}{z^s}}$$

Surely, the arguments given in this article regarding the connection to k/x distribution should apply here as well, and even more perfectly and smoothly so, hence this points to the same logarithmic proportions of LOG(1 + 1/d); yet no one would dare claiming that the integers are Benford! They are not!

In reality prime numbers do not care much about integral-powers-of-ten intervals; they merrily march forward along the integers at their own pace, totally disregarding them; nor do they have any special esteem for number 10 – an annoying member of their archrival and competitive group - the composites. They live and behave as if the number 10 does not matter much. Surely they also do not pay any attention whatsoever to the particular number systems invented by the various civilizations scattered about randomly across the universe, as they float eternally up there well above all such lowly and arbitrary local inventions. Conclusions about the primes in digital Benford's Law, either in base 10 or even employing a generic base for all positional number systems, do not interest the primes much. Yet the primes are quite enthusiastic about what was found in connection with the General Law of Relative Quantities (Kossovsky (2014), the 7th section) which assigns the primes the bin-equality proportion of 1/D ultimately at infinity, and which is number system invariant.

## Acknowledgement:

The author wishes to thank the distinguished mathematician George Andrews and Robert Hall for their helpful comments and suggestions.

March 8, 2016
Alex Ely Kossovsky
akossovsky@gmail.com